# Factoring Middle Binomial Coefficients

Gennady Eremin
ergenns@gmail.com

March 5, 2020

**Abstract**. The article describes prime intervals into the prime factorization of middle binomial coefficient. Prime factors and prime powers are distributed in layers. Each layer consists of non-repeated prime numbers which are chosen (not calculated) from the noncrossing prime intervals. Repeated factors are formed when primes are duplicated among different layers.

## 1. Introduction

**1.1. Chebyshev intervals.** This article arose thanks to the paper of Pomerance [1] that deals with the middle binomial coefficient (MBC) $B(n) = \binom{2n}{n}$. These coefficients are located in the center of the even-numbered rows of Pascal's triangle. The author notes that $\binom{2n}{n}$ is divisible by the product of all primes in the interval $(n, 2n)$. This fact was considered by Chebyshev still in 1850, and we will show that a similar *prime interval* (or the *Chebyshev interval*) is not the only one in the prime factorization of $B(n)$.

MBC's belong to the family of Catalan-like numbers, such as Catalan numbers, Motzkin numbers, hexagonal numbers, etc. The MBC is the crucial character in the definition of the Catalan number

$$C(n) = B(n)/(n+1).$$

In [2] the prime factorization of Catalan numbers is carried out using the identical Chebyshev intervals. We will apply use same techniques and methods in this paper. For the MBC the general term is

(1.1) $$B(n) = (2n)!/(n!)^2, \; n \geq 0.$$

The first coefficients are 1, 2, 6, 20, 70, 252, 924, 3432, … (sequence A000984). From (1.1) immediately follows:

(a) Prime factors of $B(n)$ are less than $2n$, i.e. these primes are contained in the open interval $(1, 2n)$, *factor-base* of $B(n)$.
(b) Any prime $p \in (n, 2n) = S_1$ divides $B(n)$ or $p \mid B(n)$. Let us call $S_1$ the *main prime interval* since it is most extensive in the factor-base.

(c) Any prime $p \in S_1$ is a non-repeated factor of $B(n)$, i.e. $p^2 \nmid B(n)$.

The next prime interval $S_2 = (\frac{1}{2}n, \frac{2}{3}n)$ is the adjacent to $S_1$. In $S_2$ all primes also divide $B(n)$, and these primes are distinct too. Between $S_1$ and $S_2$, there is the first close *factor-free zone* $Z_1 = [\frac{2}{3}n, n]$ into which all primes do not divide $B(n)$. Let us explain the above by an example.

**Example 1.1.** Prime factorization of the $1000^{th}$ MBC is the multiset which contains 217 primes and prime powers. The main prime interval (1000, 2000) "furnishes" 135 distinct prime factors to the multiset. 26 distinct primes come from the second prime interval (1000/2 = 500, 2000/3 = 666.6) additionally. Between these intervals there is the close zone [666.6, 1000] in which there are no prime factors of $B(1000)$. □

As you can see, using the Chebyshev intervals you can select a large group of distinct prime factors (not performing the usual calculations). Obviously a chain of prime intervals (separated by factor-free zones) continues with the further descent on the factor-base. As a result, prime intervals allow you to get the absolute most of prime factors.

However, a natural question arises: how to find the remaining prime factors and prime powers? And here it is logical to assume the existence of other networks of the Chebyshev intervals; otherwise it is impossible to obtain, for example, prime powers.

**1.2. Legendre layers.** For a prime $p$ and a positive integer $m$, let $v_p(m)$ denotes the largest power of $p$ dividing $m$ (see [1]). For example, $v_{11}(22) = 1$, $v_{11}(98) = 0$, and $v_{11}(23!) = 2$. This function is extended to positive rational numbers as follows:

$$v_p(a/b) = v_p(a) - v_p(b), \text{ if } v_p(a) \geq v_p(b).$$

Here well-known Legendre's formula

$$v_p(m!) = \sum_{j \geq 1} \lfloor m/p^j \rfloor.$$

According to (1.1), we have for each prime $p$

(1.2) $$v_p(B(n)) = \sum_{j \geq 1} \left( \lfloor 2n/p^j \rfloor - 2 \lfloor n/p^j \rfloor \right).$$

It is easy to show that each term in the above sum can take only the binary value 0 or 1. Obviously $\lfloor 2n/p^j \rfloor - 2 \lfloor n/p^j \rfloor = 1$ if and only if the fractional part $\{n/p^j\} \geq \frac{1}{2}$, $j \geq 1$. (See remarkable Kummer's theorem in [1].)

The equality (1.2) clearly shows the possible split of prime factors and prime powers into layers, the *Legendre layers*. Each layer contains only distinct (non-repeated) prime numbers.

**Proposition 1.1.** *Let $p$ be a prime and let $\{n/p^j\} \geq \frac{1}{2}$. Then $p$ is distributed in the jth Legendre layer of $B(n)$.*

Obviously the repeating prime factors are distributed in several layers. Pay attention to the fact that the distinct (non-repeated) prime does not necessarily fall into the first Legendre layer. For example, for $B(10^6)$ the distinct odd prime factor 11 falls into the 6[th] layer.

For the MBC with index $n$, there are only non-repeated prime factors above the border $\sqrt{2n}$. Therefore, into the interval ($\sqrt{2n}$, $2n$) there are only the distinct odd prime factors, and all these factors are distributed into the first layer (here other layers don't spread). Let's call this layer a *distinct-layer*. The first layer is the most numerous; we can say that the whole multi-layered pyramid rests on this bottom layer, so we'll start with the distinct-layer.

**Note**. In this work, we consider only odd simple factors of the MBC. The single even prime 2 is easy to calculate using (1.2). In this case for $B(n)$, the number of layers does not exceed $\log_3 n$. Legendre's formula is responsible for filling the layers; in this article, we are interested in grouping the elements of each layer into the Chebyshev intervals.

## 2. Chebyshev intervals into the distinct-layer

In this section we will consider the distribution of prime factors of $B(n)$ into the most extensive distinct-layer. The main interval $S_1 = (n, 2n)$ contains more than half of the prime factors. We have already noted the adjacent second interval $S_2 = (\frac{1}{2}n, \frac{2}{3}n)$. Let us calculate the general form of the prime intervals for the distinct-layer.

The equation (1.1) is often simplified for performing practical calculations. Let's separate the even and odd factors in the numerator of the fraction:

$$(2n)! = 2 \cdot 4 \cdot 6 \cdots 2n \times 1 \cdot 3 \cdot 5 \cdots (2n-1) = 2^n \times n! \times (2n-1)!!.$$

After the cuts it will receive (by analogy with the Catalan numbers [3])

(2.1) $$B(n) = 2^n \times (2n-1)!! \big/ n! = 2^n \times \frac{1 \cdot 3 \cdot 5 \cdots (2n-1)}{1 \cdot 2 \cdot 3 \cdots n}.$$

As a result, we get these two mathematical dependencies:

$$v_2(B(n)) = v_2(2^n) - v_2(n!) = n - \sum_{j \geq 1} \lfloor n/2^j \rfloor$$

and

(2.2) $$v_p(B(n)) = v_p((2n-1)!!) - v_p(n!)$$
$$= \tfrac{1}{2} \sum_{j \geq 1} (\lfloor (2n-1)/p^j \rfloor_{\text{odd}} + 1) - \sum_{j \geq 1} \lfloor n/p^j \rfloor, \; p > 2.$$

The calculating of the largest power of an odd prime factor for the odd double factorial and the description of the *odd floor function* $\lfloor x \rfloor_{odd}$ (the maximum odd integer not exceeding $x$) see in [2]. Since we are interested in the Chebyshev intervals, i.e. groups of adjacent prime numbers (but not individual prime factors), we will use the formula (2.1) more often.

Recall, into the main Chebyshev interval $(n, 2n)$, each prime divides $B(n)$. In the denominator of (2.1), the last odd factor $n$ may be a prime. Then $n$ is in the numerator also in a single instance. So the prime $n$ does not divide $B(n)$. The odd factor $n-1$ is also present in a single copy in the both factorials, hence the prime $n-1 \nmid B(n)$ too. The situation does not change if we will continue descending the factor-base until we get to $\tfrac{2}{3}n$.

The integer $\tfrac{2}{3}n$ is always even, i.e. $\tfrac{2}{3}n$ is not a prime. Thus any prime $\tfrac{2}{3}n < q \leq n$ is contained in a single copy into both factorials of (2.1), namely:

(a) in $(2n-1)!!$, the second copy of $3q > 3 \times \tfrac{2}{3}n = 2n$ is invalid;
(b) in $n!$, the second copy of $2q > 2 \times \tfrac{2}{3}n > n$ is invalid too.

Let $p < \tfrac{2}{3}n$ be a prime. Then we get a second copy in the numerator of (2.1) due to the composite odd factor $3p < 2n$. In the case $p > \tfrac{1}{2}n$, a symmetrical copy will not appear in the denominator with the composite number $2p > n$, so $p \mid B(n)$. As a result, we received the first (possibly not last) factor-free zone $[\tfrac{2}{3}n, n]$; in this closed set, there are no prime factors of $B(n)$.

In the numerator of (2.1) for each odd prime number $p < \tfrac{2}{3}n$, there is a guaranteed second copy due to the composite odd factor $3p < 2n$. We'll get a third copy for any odd prime $p < 2n/5$ with the additional composite factor $5p < 2n$. It is easy to see that we'll get at least $k$ copies for any odd prime $p < 2n/(2k-1) = n/(k-\tfrac{1}{2})$ with the odd factors $p, 3p, 5p, \ldots, (2k-1)p$. It is easy to see that the number of copies is more than $k$ in the case

$$2k-1 \geq p \quad \text{or} \quad k \geq (p+1)/2.$$

Thus

(2.3) $\qquad v_p((2n-1)!!) \geq k$ for any odd prime $p < n/(k-\tfrac{1}{2})$.

Obviously for $p < 2n/(2k+1) = n/(k+\tfrac{1}{2})$ we get an extra copy, i. e.

(2.4) $\qquad v_p((2n-1)!!) \geq k+1$ for any odd prime $p < n/(k+\tfrac{1}{2})$.

Using (2.3) and (2.4) easy get the following proposition.

**Proposition 2.1.** *Let $p$ be an odd prime number and let $p \in (n/(k+\tfrac{1}{2}), n/(k-\tfrac{1}{2}))$, $1 \leq k < (p+1)/2$. Then*
$$v_p((2n-1)!!) = k.$$

In Proposition 2.1, the interval bounds can't be an odd prime (only even integers).

Next, consider $n!$. For a prime $p \leq n/2$, there is a guaranteed second copy due to the composite factor $2p$. We'll get a third copy for any prime $p \leq n/3$ with the additional composite factor $3p$. Thus we'll get at least $k$ copies for a prime $p \leq n/k$ with the factors $p, 2p, 3p, \ldots, kp$. The number of copies exceeds $k$ if $k \geq p$. Thus

(2.5) $\qquad v_p(n!) \geq k$ for any odd prime $p \leq n/k$.

Obviously for $p < n/(k-1)$, the number of copies will decrease by 1, i. e.

(2.6) $\qquad v_p(n!) \geq k-1$ for any odd prime $p \leq n/(k-1)$.

Using (2.5) and (2.6) easy get another proposition.

**Proposition 2.2.** *Let $p$ be an odd prime and let $p \in (n/k, n/(k-1)]$, $1 \leq k < p$. Then*
$$v_p(n!) = k-1.$$

In case $v_p(n!) = 0$ we receive the conditional infinite interval $(n, \infty]$.

It is easy to see that in Proposition 2.1 and Proposition 2.2, the intervals for the prime $p$ partially overlap each other. The boundary points of the gaps are shown below on a straight line segment (open borders are indicated by "punctured points", and the closed upper boundary of a half-interval by a solid point). The boundaries of the interval in the double factorial are shown in red.

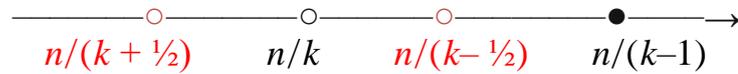
$\qquad n/(k+½) \qquad n/k \qquad n/(k-½) \qquad n/(k-1)$

To get the $k$-th Chebyshev segment of the first layer, it is enough to superimpose both sections on top of each other, i.e. intersect both intervals:

$$(n/(k+½), n/(k-½)) \cap (n/k, n/(k-1)] = (n/k, n/(k-½))$$

The following statement declares a series of Chebyshev intervals for the first layer.

**Proposition 2.3.** *Let $p$ be an odd prime number and let $k < (p+1)/2$ be a positive integer. Then*
$$v_p(B(n)) = 1 \text{ if and only if } p \in (n/k, n/(k-½)).$$

For $k=1$ we get the first (main) prime interval $S_1 = (n/1, n/(1-½))$. Further, going down the factor base, bypassing the prohibition zone, we find the second interval $S_2 = (n/2, n/(2-½))$. We can continue further descent along the base, sorting through other Chebyshev intervals, alternating with the exclusion zones, until we get to the boundary of multiplicity $\sqrt{2n}$, below which, in addition to single simple factors, squares are also possible (factors of multiplicity 2), cubes, etc.

According to Kummer's theorem, there are no repeating elements in the Legendre layer; therefore, we can continue the Chebyshev intervals to a multiple re-gion of the factor base. Since we consider only odd simple factors, the limiting value of the lower boundary of the intervals is as follows:

$$n/k > 2 \text{ or } k < n/2.$$

Now we can formulate the corresponding theorem.

**Theorem 2.1.** *For $B(n)$, the distinct-layer is formed by the Chebyshev intervals of the form*

$$S_k = (n/k, n/(k-½)), \ 1 \leq k < n/2.$$

For a given prime using Theorem 2.1, it is easy to determine it belongs to the distinct layer or not. This analysis is performed in two stages. For an odd prime $p$, the first step calculates the appropriate $k$th interval based of the lower boundary of the Chebyshev interval:

$$k = \lceil n/p \rceil.$$

In the second step, we compare $p$ and the upper bound of the $k$-th interval. In other words, $p$ falls into the $k$th interval if

$$p < n/(k-½).$$

Now we state the corresponding corollary

**Corollary 2.1.** *Let $p$ be an odd prime and let $k = \lceil n/p \rceil$. Then $p$ is included to the $k$th interval of the distinct layer of $B(n)$ if and only if $p < n/(k-½)$.*

We note that the inequality in Corollary 2.1 is strict (equality is impossible), since the integer $n/(k-½)$ is always even.

## 3. Square-layer, cube-layer, and others

We can say that the collection of the layers forms a pyramid. Each layer contains non-repeated prime factors of the MBC. The repeating primes are distributed in several layers. The distinct layer is located at the bottom of the pyramid, and it is the most numerous and it accumulates the absolute majority of distinct prime factors. The next *square-layer* accumulates the majority of primes that have the power 2. Next is the *cube-layer*, etc. The primes of the square-layer are distributed in the truncated factor-base $(1, \sqrt{2n})$. Accordingly, the primes of the cube-layer are distributed in the interval $(1, \sqrt[3]{2n})$ and so on.

In the previous section, we have built a series of distinct-layer prime intervals on the basis of (2.1). For example, the first interval $(n, 2n)$ is evident directly in (2.1). Let's convert (2.1) to simplify calculations for the square-layer. According to (1.2), the square-layer consists of those and only those primes $p$ for which

(3.1) $$\lfloor 2n/p^2 \rfloor - 2\lfloor n/p^2 \rfloor = 1 \text{ or } \{n/p^2\} \geq \tfrac{1}{2}.$$

In this regard, we can significantly reduce the number of factors in both factorials of (2.1). It is easy to see that in the numerator the odd factors can be limited to the value of $\sqrt{2n}$ since we are working only with squares. Indeed, in (3.1) the operand $\lfloor 2n/p^2 \rfloor$ is zero for any $p > \sqrt{2n}$. Accordingly, the factors in the denominator can be restricted to $\sqrt{n}$. The conversion in (2.1) is therefore as follows:

$$(2n-1)!!\big/n! \Rightarrow (\lfloor\sqrt{2n}\rfloor_{odd})!!\big/\lfloor\sqrt{n}\rfloor!$$
(3.2)
$$= 1\cdot 3\cdot 5\cdots\lfloor\sqrt{2n}\rfloor_{odd}\big/1\cdot 2\cdot 3\cdots\lfloor\sqrt{n}\rfloor.$$

Here $\lfloor\sqrt{2n}\rfloor_{odd}$ is the maximum odd integer less then $\sqrt{2n}$.

In (3.2) everyone can see the interval $(\sqrt{n}, \sqrt{2n})$ in which each prime belongs to the square-layer. This is the longest similar interval, i.e. we are dealing with the main Chebyshev interval of the square-layer. Obviously the main Chebyshev interval for an arbitrary $j$th layer is determined as follows

(3.3) $$S_1^{(j)} = (\sqrt[j]{n}, \sqrt[j]{2n}), \; j \geq 1.$$

All prime numbers into interval (3.3) belong to the $j$th layer of $B(n)$. Next, we will analyse by analogy with the distinct-layer.

In the denominator of (3.2), the last factor $q = \lfloor n^{1/2} \rfloor$ may be an odd prime number. Then $q$ is in the numerator and also in a single instance. Therefore $q$ does not belong to the square-layer. Note $q$ can divide $B(n)$ if it belongs to the distinct-layer (the cube-layer and other layers do not extend here).

Now let's look at the smaller factors of the form $q = \lfloor (n-i)^{1/2} \rfloor$, $i = 1, 2, \ldots$ For the initial values of the variable $i$, a possible prime $q$ is still present in a single copy in the both factorials. The situation will not change until $i \leq \tfrac{1}{3}n$ (any integer $(\tfrac{2}{3}n)^{1/2}$ always even and therefore composite). Thus each odd prime $(\tfrac{2}{3}n)^{1/2} < q \leq n^{1/2}$ is contained in the single copy into the both factorials of (3.2), so it does not belong to the square-layer. Indeed, on the one hand in the numerator of (3.2), the next copy with $3q^2 > 3\times\tfrac{2}{3}n = 2n$ is invalid; and on the other hand in the denominator, the second copy with $2q^2 > 2\times\tfrac{2}{3}n = n + \tfrac{1}{3}n$ is invalid too.

As a result, we received the first factor-free zone $Z_1^{(2)} = [(\tfrac{2}{3}n)^{1/2}, n^{1/2}]$ for $B(n)$. Each prime in this close zone does not belong to the square-layer, but inside $Z_1^{(2)}$ can be elements of the distinct-layer. It is easy to show that the first factor-free zone for an arbitrary $j$th layer is determined as follows:

$$Z_1^{(j)} = [(\tfrac{2}{3}n)^{1/j}, n^{1/j}], \; j \geq 1.$$

Now we consider the second prime interval of the square-layer and its components. Let's look at the odd prime factors of the form $(\frac{1}{2}n)^{1/2} < p < (\frac{2}{3}n)^{1/2}$. In the numerator of (3.2) we get the second copy due to the composite factor $3p^2 < 3 \times \frac{2}{3}n = 2n$; at the same time in the denominator, the second copy with the factor $2p^2 > n$ is invalid. So we got the second Chebyshev interval of the square-layer

$$S_2^{(2)} = ((\tfrac{1}{2}n)^{1/2}, (\tfrac{2}{3}n)^{1/2}) \quad \text{or for any } j\text{th layer} \quad S_2^{(j)} = ((\tfrac{1}{2}n)^{1/j}, (\tfrac{2}{3}n)^{1/j}).$$

Write the dependencies for the first and second intervals in the following form:

$$S_1^{(2)} = ((n/1)^{1/2}, (n/(1-\tfrac{1}{2}))^{1/2}) \quad \text{and} \quad S_2^{(2)} = ((n/2)^{1/2}, (n/(2-\tfrac{1}{2}))^{1/2})).$$

Obviously

$$S_3^{(2)} = ((n/3)^{1/2}, (n/(3-\tfrac{1}{2}))^{1/2}) \quad \text{or for any } k\text{th interval} \quad S_k^{(2)} = ((n/k)^{1/2}, (n/(k-\tfrac{1}{2}))^{1/2})).$$

The trend is clear and we can write a formula for the $k$th interval of any $j$th layer:

(3.4) $$S_k^{(j)} = ((n/k)^{1/j}, (n/(k-\tfrac{1}{2}))^{1/j})).$$

For the factorization of $B(n)$, we don't need to consider all the layers the number of which can reach $\log_3 n$. Usually, it is sufficient to select the Chebyshev intervals of some lower layers of the pyramid. For example, the prime factorization of $B(10^6)$ consists of 101538 primes and prime powers, of which 101384 fall into the distinct-layer including 70435 primes in the main interval $(10^6, 2 \cdot 10^6)$. Using intervals of the square-layer 120 primes are selected. The remaining 34 primes and prime powers are easily calculated from (1.2) or (2.2).

## 4. Software service

In conclusion, let us consider a software system for the implementation of test calculation in the real time. The following is a list of several programs with brief description of their functions.

(1) A small additional program can produce a list of primes in a certain range; the start of the range cannot exceed $10^{10}$.

(2) The following program executes the interval factorization of the MBC for a given index $n$. Prime factors and power primes are issued in ascending order starting with 2 and ending primes of the main interval $(n, 2n)$. Older factors are grouped into prime intervals as possible. Each factor is accompanied (in parentheses) by a list of layers on which this prime is distributed. The first layer for distinct (non-repeated) factors is not indicated. For example, this applies to all primes factors exceeding the power border $\sqrt{2n}$.

(3) The user can obtain the layer-splitting of the prime factors of any MBC. Each layer contains only unique (non-repeated) primes. On the printout you can see how repetitive prime factors are divided into different layers. The layers are processed starting from the top and finishing the bottom distinct-layer. Separately it is allocated empty layers that contain no primes.

(4) It is possible to calculate a separate layer. You can display separately square-layer and cube-layer. Working with an individual layer, it is possible to obtain the information for the MBC with very large indexes (billion and more).

**Acknowledgements.** I wish to thank Igor Pak (Mathematics Department at UCLA, USA) for the detailed and comprehensive historical review of Catalan-like numbers, and thanks for a selection of recent works on the subject. I am grateful to Bruce Sagan (Michigan State University, USA) who was interested in the prime intervals into the prime factorization of Catalan numbers and MBC's, and who approved of the developed software service.

Gzhel State University, Moscow, 140155, Russia
http://www.en.art-gzhel.ru/